\documentclass[12pt]{amsart}
\usepackage{amssymb}
\usepackage{amscd}
\usepackage{color}
\usepackage{url}
\usepackage{hyperref}
\usepackage{xcolor}

\numberwithin{equation}{section}

\def\re{\mathop{\rm Re}\nolimits}

\theoremstyle{plain}
\newtheorem{theorem}{Theorem}[section]

\newtheorem{thm*}{Theorem}

\newtheorem{cor}[theorem]{Corollary}

\theoremstyle{definition}

\newtheorem{ex}[theorem]{Example}
\theoremstyle{remark}
\newtheorem{rem}[theorem]{Remark}

\def\CC{\mathbb C}

\def\NN{\mathbb N}
\def\RR{\mathbb R}
\def\SS{\mathbb S}

\def\beginpf{\begin{proof}}
\def\endpf{\end{proof}}
\def\beq{\begin{equation}}
\def\eeq{\end{equation}}

\newcommand{\xdownarrow}[1]{%
  {\left\downarrow\vbox to #1{}\right.\kern-\nulldelimiterspace}
}

\begin{document}

\title[Weighted composition operators and monomial operators]{Boundedness and compactness of weighted composition operators and monomial operators}
\author{I. Chalendar}
\address{Isabelle Chalendar,  Universit\'e Gustave Eiffel, LAMA, (UMR 8050), UPEM, UPEC, CNRS, F-77454, Marne-la-Vallée (France)}
\email{isabelle.chalendar@univ-eiffel.fr}

\author{J.R. Partington}
\address{Jonathan R. Partington, School of Mathematics, University of Leeds, Leeds LS2 9JT, UK}
 \email{J.R.Partington@leeds.ac.uk}

\maketitle

{\bf Abstract:}
This paper characterises the boundedness and compactness of Agler--McCarthy monomial
operators
by reducing them to weighted composition operators and 
deriving explicit Carleson measure criteria on the
half-plane. The results are illustrated by examples.\\


{\bf MSC class:} 30H10, 47B33

\section{Introduction}

Motivated by questions of polynomial approximation,
Agler and McCarthy \cite{AM} defined a {\em monomial operator\/}
as a bounded  linear
operator $T: L^2(0,1) \to L^2(0,1)$ mapping monomials to monomials;
that is, 
such that for each $n \in \NN$ we have (with a slight abuse of notation) $T(x^n)=c_n x^{p_n}$ for some $p_n \in \CC$ with $\re p_n > -1/2$.
Such an operator is {\em flat\/} if there is a constant $b$ with $p_n=n+b$
for each $n$. We shall also say that it has {\em affine index\/} if $p_n=an+b$ for
real constants $a,b \ge 0$.

Well-known examples of monomial operators include the Volterra operator
$V$ with 
\[
Vf(x)=\int_0^x f(t) \, dt, \qquad Vx^n=\frac{x^{n+1}}{n+1},
\]
and the Hardy operator $H$ (also known as the continuous Ces\`aro operator \cite{BHS}) with
\[
Hf(x)=\frac{1}{x} \int_0^x f(t) \, dt, \qquad Hx^n=\frac{x^{n}}{n+1},
\]
which are clearly flat operators. Agler and McCarthy gave conditions for flat operators
to be bounded, and showed that they were never compact except in trivial cases,
but the general case was left open.

Their work was based on a transformation to the Hardy space of the
half-plane $\SS=\{s \in \CC: \re s > -1/2 \}$ and showed that the adjoint of
a bounded monomial operator can be expressed as an operator
on $H^2(\SS)$. By means of a closely-related transformation we shall work on the
Hardy space $H^2(\CC_+)$ on the standard right half-plane, and derive results
in the theory of weighted composition operators, which will be applied to characterise
bounded and compact operators. The results are particularly simple
for flat and affine-index operators, as will be illustrated by examples.

\section{Boundedness and Compactness}

The approach here is very similar to Agler--McCarthy \cite{AM}, but adapted for own purposes.\\

We transfer the operator $T$ to $L^2(0,\infty)$. Writing $x=e^{-t}$  we have 
\[
\int_0^1 |f(x)|^2 \, dx = \int_0^\infty |f(e^{-t})|^2 e^{-t} \, dt,
\]
giving an isometry by $Jf(t)= e^{-t/2} f(e^{-t})$. So $J^{-1}g(x) = \frac{1}{\sqrt x} g(-\log x)$.
The unitarily equivalent operator $\tilde T$ on $L^2(0,\infty)$ maps $e^{-nt-t/2}$ to
$c_n e^{-p_nt-t/2}$.

Then by taking Laplace transforms and using the
Paley--Wiener theorem we have a unitarily equivalent  operator $T_0$ on $H^2(\CC_+)$ which
maps $1/(s+n+1/2)$ to $c_n /(s+p_n+1/2)$.

Since the reproducing kernel for $H^2(\CC_+)$ is $k_w(s)=1/(s+\overline w) $ we have
\[
T_0^* G(n+1/2) = \overline{c_n} G(p_n+1/2).
\]
The Blaschke sequences (zero sequences) $(z_n)$ for $H^2(\CC_+)$
are those with $\displaystyle \sum_n \frac{\re z_n}{1+|z_n|^2}<\infty$
\cite{duren}, and so
the sequence $(n+1/2)$ is not a Blaschke sequence for the
Hardy space of the right half-plane. Thus we have a unique formula for $T_0^*$,
namely
 the weighted composition operator $W_{h,\phi}$
with
\beq\label{eq:phih}
\phi(n+1/2)=p_n+1/2 \qquad \hbox{and} \qquad h(n+1/2)=\overline{c_n},
\eeq
since $(T_0^*-W_{h,\phi})(G)=0$ for all $G \in H^2(\CC_+)$.

For affine-index operators with $p_n=an+b$ we have
$T_0^* G=h( G\circ \phi)$ with $h(n+1/2)=\overline c_n$ and $\phi(s)=as+b$.
This is consistent with Agler--McCarthy \cite{AM} except that we use the ``usual'' right half-plane
not the shifted version.

\begin{rem}
The classical M\"untz--Sz\'asz theorem says that the functions $\{x^n: n \in S \}$
span a dense subspace of $L^2(0,1)$ if and only if $ \displaystyle \sum_{n\in S} \frac{1}{n+1} = \infty$. Thus under this condition a bounded monomial operator is completely defined by
$\{T x^n: n \in S\}$.
We see this from the above discussion, since such an $S$ does not constitute
a Blaschke sequence, and the interpolation problem has a unique solution.
\end{rem}

Boundedness of the weighted composition operator $W_{h,\phi}: f \mapsto h (f \circ \phi)$ can be characterised using a Carleson measure criterion by proving a result similar to those in \cite{CHD,KP}.
Namely, for Borel subsets $E \subseteq \overline{\CC_+}$ define a measure by
\begin{equation}\label{eq:defmu}
\mu(E)= \int_{\phi^{-1}(E) \cap i\RR} |h(s)|^2 \, |ds|.
\end{equation}

\begin{theorem}
The weighted composition operator $W_{h,\phi}$ is bounded
 if and only if  $\mu$ is a Carleson measure; that is,
$\sup_I \mu(Q_I)/|I| < \infty$, where  $Q_I=[0, |I|] \times I$
is a Carleson square based on an interval $I \subset i\RR$.
Hence the   monomial operator $T$ is bounded if and only if the same condition holds
for the associated weighted composition operator.
\end{theorem}

\beginpf
In \cite{CHD} the analogous result is proved for weighted composition
operators on the disc
and in \cite{KP} a similar result is derived for the Bergman space.
For the half-plane an identical proof works with obvious modifications, as follows:

We claim that for $g \ge 0$ Borel measurable we have
\begin{equation}\label{eq:normid}
\int_{\overline{\CC_+}} g \, d\mu = \int_{i\RR} |h|^2 (g \circ \phi) |ds|.
\end{equation}
As in \cite{CHD} we can easily verify this result for simple functions, and then 
for arbitrary $g$ take an increasing sequence of positive simple functions $(g_n)$
with $g_n \to g$ pointwise in $\overline{\CC_+}$.

From \eqref{eq:normid} we deduce that $\|W_{h,\phi}f\|= \|f\|_{L^2(\overline{\CC_+},\mu)}$.
Now, the Carleson measure condition that for some fixed $C>0$
\[
\int_{\overline{\CC_+}} |f|^2 \, d\mu \le C \|f\|^2 \qquad (f \in H^2(\CC_+))
\]
 is equivalent to the condition that $W_{h,\phi}$ is bounded.
 
In the context of a monomial operator $T$ the adjoint operator
is unitarily equivalent to 
$W_{h,\phi}$ and the result follows.
 \endpf
 
 \begin{rem}
 It is interesting to consider the case $h=1$, an unweighted composition
 operator $C_\phi$ on the right half-plane. We see that 
 $C_\phi$ is bounded if and only if there is a constant $K>0$ such that
 $\mu(Q_I)= m(\phi^{-1}(E) \cap i\RR) \le K|I|$
 for all Carleson squares (here $m$ is Lebesgue measure). This may be compared with the
 alternative condition given by Elliott and Jury \cite{EJ}; namely, that there is an 
 angular derivative  at $\infty$, the non-tangential limit  
 $\phi'(\infty)=\lim_{z \to \infty\ {\rm n.t.}}z/\phi(z)$ satisfying $0 < \phi'(\infty)< \infty$.
 \end{rem}

\begin{cor}
Let $T$ be an  affine-index monomial operator with $p_n=an+b$ ($a,b>0$).
Then $T$ is bounded if and only if the associated function $h$ satisfies
\beq\label{eq:cclin}
\sup_{t \in \RR} \int_{t}^{t+1} |h(iy)|^2   \,  dy < \infty.
\eeq
\end{cor}

\beginpf
By the simple structure of
$\phi$, we have
$\mu$ concentrated on $b+i\RR$ and the Carleson condition is
\[
\int_{u/a}^{v/a} |h(iy)|^2 \, dy \le C (v-u)
\]
but only for $v-u \ge b$; that is, equivalently,
\eqref{eq:cclin}.
\endpf

\begin{ex}{\rm 
We can have $c_n \to 0$ and the monomial operator not even bounded.
For example, take $c_n=(n+1/2)e^{-(n+1/2)}$, for which the associated function
$h$ can only be $se^{-s}$, which does not satisfy the Carleson criterion.
This fact may also
be proved using the fact that the monomials do not form an unconditional basis, but this is more
explicit.}
\end{ex}

We know from \cite{AM} that there are no compact monomial operators
of the form $x_n \mapsto c_n x^n$, except the zero operator, and we also
know that
there are no compact composition operators on the
half-plane \cite{matache}.\\

With $T(x^n)=c_n x^{p_n}$, if $T$ is bounded then the sequence
\begin{equation}\label{eq:txn}
\|Tx^n\|/\|x^n\|=|c_n| \sqrt{\frac{2n+1}{2\re p_n+1}}
\end{equation}
is bounded. So for affine-index  operators  if $T$ is bounded we clearly have $(c_n)$ bounded.

If $T$ is compact, then observe first that
 the normalized functions $e_n : x \mapsto \sqrt{2n+1} \, x^n$ tend weakly to zero.
For they form a bounded sequence  tending weakly to zero on the dense subspace consisting of polynomials, since
\[
\langle e_m, e_n \rangle = \frac{\sqrt{2m+1}\sqrt{2n+1}}{m+n+1} \to 0 \qquad \hbox{as} \quad n \to \infty.
\]
Hence compactness implies that $\|Te_n\| \to 0$ and this is 
given by \eqref{eq:txn}. That is, for compact affine-index operators, $c_n \to 0$.
\\

Using the completely continuous characterization of compactness in
reflexive Banach spaces,
namely, that if $f_n \to 0$ weakly then $\|W_{h,\phi}f_n \| \to 0$,
we see that 
a necessary and sufficient condition on compactness
will be  compactness of the mapping 
$J: H^2(\CC_+) \to L^2(\overline{\CC_+},\mu)$. This is characterised
 for the disc in \cite{CHD}.

On the disc this is the familiar  vanishing
Carleson measure criterion for $\mu$.
On the half-plane the condition that 
\beq\label{eq:vcm}
\mu(Q_I)/|I| \to 0 \qquad \hbox{as} \qquad |I| \to 0,
\eeq where $Q_I=[0, |I|] \times I$
is a Carleson square based on an interval $I \subset i\RR$,  is not  
sufficient, as we now show.\\

Consider the mapping $x^n \mapsto x^{n+1}$, which is bounded but
not compact. We have $h(s)=1$ and $\phi(s)=s+1$.
Now according to \eqref{eq:defmu} we have $\mu(E)=m(E \cap (1+i\RR))$, for $E \subset \CC_+$,
and thus
\[
\|C_\phi f\|^2 = \int_{\overline{\CC_+}} |f|^2 \, d\mu = \int_{1+i\RR} |f|^2 \, |ds|.
\]
 The measure $\mu$ of a square $Q_I$ is $0$ if $|I|<1$ and so   this
 $\mu$ satisfies \eqref{eq:vcm}, but nonetheless we do not have compactness.
 \\

 In fact the 
correct test for compactness for embeddings from $H^2(\CC_+)$ into $L^2(\overline{\CC_+},\mu)$  
is given in \cite{PPW} in the more general context of Zen spaces, namely, it can be written
\beq\label{eq:truevc}
\sup_{I: C(I) \in S_r} \mu(Q_I)/|I| \to 0 \qquad \hbox{as} \quad r \to 0,
\eeq
where $C(I)$ is the centre of a Carleson square $Q_I$ and for $0<r<1$ we define
\[
S_r = \{z \in \CC_+: 0<\re z<r \hbox{ or } |z|>1/r\},
\]
which is
the union of the complement of
a large semi-circle and a narrow strip near the axis.

Equivalently, for each $\epsilon>0$ there is a compact subset $K \subset \CC_+$ such that
$\mu(Q_I)/|I| < \epsilon$ for every $Q_I$ with centre lying outside $K$.\\

We therefore have the following result.
\begin{theorem}
The monomial operator $T$ is compact if and only if the
associated measure $\mu$ satisfies the condition \eqref{eq:truevc}.
\end{theorem}

This can be made more explicit for operators with affine index.

\begin{cor}\label{cor:cpct}
Let $T$ be an  affine-index monomial operator with $p_n=an+b$ ($a,b>0$).
Then $T$ is compact if and only if the associated function $h$ satisfies the condition that for each
$\epsilon>0$ there is an $M>0$ such that
\beq\label{eq:vanL}
\int_{(t-L)/a}^{(t+L)/a} |h(iy)|^2   \,  dy  < \epsilon L
\eeq
for all  $|t| \ge M$ and $L \ge b/2$. 
This holds whenever $h \in L^2(i\RR)$.
\end{cor}

\beginpf
Equation \eqref{eq:vanL} is a direct
rewriting of the  vanishing condition   \eqref{eq:truevc}, using the fact that
a square $[0,2L] \times [t-L,t+L]$ with centre $L+it$ meets the line $\{s \in \CC: 
\re s=b\}$ if and only if $L \ge b/2$. 
Now if $h \in L^2(i\RR)$, the inequality \eqref{eq:vanL} holds for all $t$ if
$L > L_0:= \|h\|_2^2/\epsilon$. For $b/2 \le L \le L_0$, we may make
the integral less than $\epsilon b/2$ by choosing $|t|$ sufficiently large. 
\endpf

\section{Examples}

For the Volterra operator,  we 
have from \eqref{eq:phih} that $h(s)=(s+1/2)^{-1}$ and $\phi(s)=s+1$.
 Now we can apply Corollary \ref{cor:cpct} and since
\[
\int_{-\infty}^{\infty} \frac{1}{|1/2+iy|^2 } \, dy < \infty
\]
we can deduce that the operator is compact, as is well known.\\

For the bounded but non-compact Hardy operator we
have $h(s)=(s+1/2)^{-1}$ and $\phi(s)=s$.
Corollary \ref{cor:cpct} does not apply since $a=0$, and we must consider all Carleson squares,
since they all meet the imaginary axis. We   have
\[
\int_u^v  |h(iy)|^2 \, dy = O(v-u), \qquad \hbox{but not} \qquad o(v-u)
\]
as $v-u \to 0$, and so the operator is bounded but not compact, as is again well known.\\

For some non-flat examples consider the following operators defined on $L^2(0,1)$:
\begin{eqnarray*}
T_1 f(x) &=& f(x^2), \\
T_2 f(x) &=& \int_0^x f(t^2) \, dt,\\
T_3 f(x) &=& \int_0^x t f(t^2) \, dt.
\end{eqnarray*}
We have
$T_1 (x^n)=x^{2n}$, $T_2 (x^n)=\dfrac{x^{2n+1}}{2n+1}$,
$T_3 (x^n)= \dfrac{x^{2n+2}}{2n+2}$.
These operators are not flat, but they are affine-index.

It is elementary to see that $T_1$ is unbounded, since it sends $x^{-1/3}$ to $x^{-2/3}$.
We shall see that $T_2$ is unbounded, but $T_3$ is bounded (even compact).\\

For $T_1$ we have $\phi(n+1/2)=2n+1/2$, giving $\phi(s)=2s-1/2$, which is not a
self-map of $\CC_+$.\\

For $T_2$ 
take $f(x)=1/(\sqrt t \log t)$ so that $f \in L^2(0,1)$.
Then 
\[
T_2 f(x)= \int_0^x \frac{1}{t (2 \log t)},
\]
and since this integral diverges we see that $T_2$ is unbounded.
Now we have $\phi(n+1/2)=2n+3/2$, giving $\phi(s)=2s+1/2$. Also
$h(n+1/2)=1/(2n+1)$, giving $h(s)=1/(2s)$.  

We may also use the Carleson test
so that  $\phi^{-1}(E)\cap i\RR=\phi^{-1}(E  \cap (1/2+i \RR))$. However since $|h|^2$ is not locally
integrable on the imaginary axis we clearly do not get a Carleson measure.\\

For $T_3$ we have $\phi(n+1/2)=2n+5/2$, and $\phi(s)=2s+3/2$.
Then $h(n+1/2)=1/(2n+2)$, giving $h(s)=1/(2s+1)$.
It is now easy to see that $\mu$ is a Carleson measure, 
concentrated on the line $\{s: \re s=3/2\}$,
and the operator is bounded.
Indeed,   since
$h \in L^2(i\RR)$ we see from Corollary \ref{cor:cpct} that $T_3$ is indeed compact.


\begin{thebibliography}{99}

\bibitem{AM}
J. Agler and J.E. McCarthy,  Monomial operators. 
{\em Acta Sci. Math. (Szeged)\/} 88 (2022), no. 1--2, 371--381. 

\bibitem{BHS}
A. Brown, P.R. Halmos and A.L. Shields,   Ces\`aro operators. {\em Acta Sci. Math. (Szeged)\/} 26 (1965), 125--137.

\bibitem{CHD}
M.D. Contreras and A.G. Hern\'andez-D\'\i az,   Weighted composition operators on Hardy spaces. 
{\em J. Math. Anal. Appl.} 263 (2001), no. 1, 224--233.

\bibitem{duren}
P.L. Duren,  
{\em Theory of $H^p$ spaces}.
Academic Press, New York-London, 1970.

\bibitem{EJ}
S. Elliott and M.T. Jury,   Composition operators on Hardy spaces of a half-plane. 
{\em Bull. Lond. Math. Soc.} 44 (2012), no. 3, 489--495.

\bibitem{KP}
R. Kumar and J.R.  Partington,  Weighted composition operators on Hardy and Bergman spaces. 
{\em Recent advances in operator theory, operator algebras, and their applications}, 157--167, Oper. Theory Adv. Appl., 153, Birkh\"auser, Basel, 2005. 

\bibitem{matache}
V. Matache,   Composition operators on Hardy spaces of a half-plane. 
{\em Proc. Amer. Math. Soc.} 127 (1999), no. 5, 1483--1491.

\bibitem{PPW}
C. Pang, A. Per\"al\"a and M. Wang,   Embedding theorems and area operators on Bergman spaces with doubling measure. 
{\em Complex Anal. Oper. Theory\/} 15 (2021), no. 3, Paper No. 42, 24 pp.


\end{thebibliography}
\end{document}